\documentclass[10pt]{article}

\usepackage{a4wide}
\usepackage{amssymb}
\usepackage{amsfonts}
\usepackage{amsmath}
\input xy
\xyoption{arrow} \xyoption{matrix}

\date{}

\newtheorem{proposition}{Proposition}[section]
\newtheorem{theorem}[proposition]{Theorem}

\def\der{\partial }

\def\nFM0{{\nu }_{F,M_0}}
\def\nFN0{{\nu }_{F,N_0}}
\def\nGN0{{\nu }_{G,N_0}}

\def\N0{ {\bf N}_0 }

\def\v{\varphi}
\def\ra{\rightarrow}

\def\Xpm{X^{\pm }}

\def\s{\sigma}

\def\l1{{\lambda}_1}

\def\a{\alpha}
\def\a0{ {\alpha }_0}
\def\a1{ {\alpha }_1}

\def\l{\lambda}
\def\o{\omega}

\def\nFGM0{{\nu }_{F,G,M_0}}

\def\nFN0{{\nu}_{F,N_0}}

\def\sm{{\sigma}^m}

\def\sm1{{\sigma}^{-1}}

\def\smtp1{{\sigma}^{-t+1}}

\def\o{\omega }
\def\S1{S^{-1}}

\def\Xpm1{X^{\pm 1}_1}

\def\sPM1{{\sigma }^{\pm 1}}
\def\sMP1{{\sigma }^{\mp 1 }}

\def\d{\delta}

\def\di{{\rm d.ind}}

\def\L{\Lambda}

\def\CA{{\cal A}}

\def\Ytm1{Y^{t-1}}
\def\Yim1{Y^{i-1}}

\def\CG{{\cal G}}
\def\CH{{\cal H}}

\def\Aut{{\rm Aut}}

\def\Der{{\rm Der }}
\def\ad{{\rm ad }}
\def\dim{{\rm dim }}

\def\CJ{ {\cal J}}

\def\SL2Z{ {\rm SL}_2({\bf Z}) }

\def\Gp1{ G^{1 , 1 } }
\def\P11{ P^{-1 , 1 } }
\def\Pp1{ P^{1 , 1 } }

\def\nCLsr{{}^\nu\kern-2pt {\cal L}^{\sigma , \rho  }}
\def\nP{{}^\nu \kern-2pt P}
\def\nL{{}^\nu\kern-2pt L}
\def\nLL{{}^\nu\kern-2pt \Lambda}
\def\nPsr{{}^\nu\kern-2pt P^{\sigma , \rho  }}
\def\nLsr{{}^\nu\kern-2pt L^{\sigma , \rho  }}
\def\nuCL{{}^\nu\kern-2pt  {\cal L}}
\def\nCLsr{{}^\nu\kern-2pt {\cal L}^{\sigma , \rho  }}
\def\nCL1m{{}^\nu\kern-2pt {\cal L}^{-1 , 1  }}
\def\x1nu{x^\frac{1}{\nu}}
\def\xm1nu{x^{-\frac{1}{\nu}}}

\def\ra{\rightarrow }

\def\CB{{\cal B}}

\def\CT{{\cal T}}

\def\CH{ {\cal H}}

\def\nAM0{{\nu }_{{\cal A},M_0}}
\def\nAN0{{\nu }_{{\cal A},N_0}}

\def\Der{ {\rm Der }}
\def\CJ{ {\cal J }}

\def\det{ {\rm det }}
\def\ad{ {\rm ad }}

\def\ga{\mathfrak{a}}
\def\gb{\mathfrak{b}}

\def\SL{{\rm SL}}

\def\di!{\frac{\der^i}{i!}}
\def\dik!{\frac{\der^k_i}{k!}}

\def\N{\mathbb{N}}

\def\0{\overline{0}}
\def\1{\overline{1}}

\def\Ln1{\L_{n,\overline{1}}}

\def\a1{a_{\overline{1}}}

\def\S{\Sigma}

\def\vn1{\overrightarrow{n-1}}

\def\mJ{\mathbb{J}}
\def\mI{\mathbb{I}}

\def\mF{\mathbb{F}}

\def\mT{\mathbb{T}}

\def\mE{\mathbb{E}}

\def\K1{{\rm K}_1}

\def\hmI1{\widehat{\mI_1}}
\def\tmI1{\widetilde{\mI_1}}
\def\tmJ1{\widetilde{\mJ_1}}
\def\hB1{\widehat{B_1}}
\def\hCB1{\widehat{\CB_1}}

\def\ggu{\mathfrak{u}}

\def\UAut{{\rm UAut}}

\def\mJ{\mathbb{J}}

\begin{document}

\author{V. V. \  Bavula   
}

\title{Every monomorphism of the Lie algebra of unitriangular polynomial derivations is an automorphism}

\maketitle

\begin{abstract}

We prove that every monomorphism of the Lie algebra $\ggu_n$  of unitriangular  derivations  of the polynomial algebra $P_n=K[x_1, \ldots , x_n]$ is an automorphism.

$\noindent $



\end{abstract}


\section{Introduction}

Throughout, 
$K$ is a
field of characteristic zero and  $K^*$ is its group of units;
$P_n:= K[x_1, \ldots , x_n]=\bigoplus_{\alpha \in \N^n}
Kx^{\alpha}$ is a polynomial algebra over $K$ where
$x^{\alpha}:=x_1^{\alpha_1}\cdots x_n^{\alpha_n}$;
$\der_1:=\frac{\der}{\der x_1}, \ldots , \der_n:=\frac{\der}{\der
x_n}$ are the partial derivatives ($K$-linear derivations) of
$P_n$;  $\Aut_K(P_n)$ is the group of automorphisms of the polynomial algebra $P_n$; $\Der_K(P_n) =\bigoplus_{i=1}^nP_n\der_i$ is the Lie
algebra of $K$-derivations of $P_n$; $A_n:= K \langle x_1, \ldots
, x_n , \der_1, \ldots , \der_n\rangle  =\bigoplus_{\alpha , \beta
\in \N^n} Kx^\alpha \der^\beta$  is  the $n$'th {\em Weyl
algebra}; for each natural number $n\geq 2$,
$$\ggu_n :=
K\der_1+P_1\der_2+\cdots +P_{n-1}\der_n$$ is the  {\em Lie algebra of unitriangular polynomial derivations} (it is a Lie subalgebra of the Lie algebra $\Der_K(P_n)$) and $G_n:= \Aut_K(\ggu_n)$ is its group of automorphisms;
$\d_1:=\ad (\der_1), \ldots , \d_n:=\ad (\der_n)$ are the inner derivations of the Lie algebra $\ggu_n$ determined by the elements $\der_1, \ldots , \der_n$ (where $\ad (a)(b):=[a,b]$).

The aim of the paper is to prove the following theorem.

\begin{theorem}\label{10Mar12}
Every monomorphism of the Lie algebra $\ggu_n$ is an automorphism.
\end{theorem}

{\it Remark}. Not every epimorphism of the Lie algebra $\ggu_n$ is an automorphism. Moreover, there are countably many distinct ideals $\{ I_{i\o^{n-1}}\, | \, i\geq 0\}$ such that $$I_0=\{0\}\subset I_{\o^{n-1}}\subset I_{2\o^{n-1}}\subset \cdots \subset I_{i\o^{n-1}}\subset \cdots$$ and the Lie algebras $ \ggu_n/I_{i\o^{n-1}}$ and $\ggu_n$ are isomorphic (Theorem 5.1.(1), \cite{Bav-Lie-Un-GEN}).

Theorem \ref{10Mar12} has bearing of  the Conjecture of Dixmier  \cite{Dix}
 for the Weyl algebra $A_n$ over a field of characteristic zero that claims: {\em every homomorphism of the Weyl algebra is an automorphism}. The Weyl algebra $A_n$ is a simple algebra, so every algebra endomorphism of $A_n$  is a monomorphism. This conjecture is open since 1968 for all $n\geq 1$. It is stably equivalent to the Jacobian  Conjecture for the polynomial algebras as was shown by  Tsuchimoto
\cite{Tsuchi05}, Belov-Kanel and Kontsevich \cite{Bel-Kon05JCDP},
(see also  \cite{JC-DP} for a short proof). The Jacobian Conjecture claims that {\em certain} monomorphisms of the polynomial algebra $P_n$ are isomorphisms: {\em Every algebra endomorphism $\s $ of the polynomial algebra $P_n$ such that $\CJ (\s ):= \det (\frac{\der \s (x_i)}{\der x_j})\in K^*$ is an automorphism.} The condition that $\CJ (\s )\in K^*$ implies that that the endomorphism $\s$ is a monomorphism.

An analogue of the Conjecture of Dixmier is true for the algebra $\mI_1:= K\langle x, \frac{d}{dx}, \int \rangle$ of polynomial
 integro-differential operators.

  \begin{theorem}\label{11Oct10}
{\rm (Theorem 1.1, \cite{Bav-cdixintdif})} Each algebra endomorphism of $\mI_1$ is an automorphism.
\end{theorem}

In contrast to the Weyl algebra $A_1=K\langle x, \frac{d}{dx} \rangle$, the algebra of polynomial
 differential operators, the algebra $\mI_1$ is neither a left/right Noetherian algebra nor a simple algebra. The left localizations, $A_{1,\der}$ and $\mI_{1, \der}$, of the algebras  $A_1$ and $\mI_1$ at the powers of the element $\der= \frac{d}{dx}$ are isomorphic. For the simple  algebra $A_{1,\der} \simeq \mI_{1, \der}$, there are algebra endomorphisms that are not automorphisms \cite{Bav-cdixintdif}.

Before giving the proof of Theorem \ref{10Mar12}, let us recall several results that are used in the proof.

{\bf The derived series for the Lie algebra $\ggu_n$}.
 Let $\CG$ be a Lie algebra over the field $K$ and $\ga$, $\gb$ be
its ideals. The {\em commutant}  $[\ga , \gb ]$ of the ideals
$\ga$ and $\gb$ is the linear span in $\CG$ of all the elements
$[a,b]$ where $a\in \ga$ and $b\in \gb$. Let $\CG_{(0)}:=\CG$,
 $\CG_{(1)}:= [\CG , \CG ]$  and $ \CG_{(i)}:=[\CG_{(i-1)}, \CG_{(i-1)}]$ for $i\geq 2$. The descending series of ideals of the Lie algebra $\CG$,
 $$ \CG_{(0)}=\CG \supseteq   \CG_{(1)}\supseteq \cdots \supseteq \CG_{(i)}\supseteq
   \CG_{(i+1)}\supseteq \cdots $$
  is called the {\em derived  series}  for the
Lie algebra $\CG$. The Lie algebra $\ggu_n$
admits the finite strictly descending chain of ideals
\begin{equation}\label{seruni}
\ggu_{n,1}:=\ggu_n\supset \ggu_{n,2}\supset \cdots \supset
\ggu_{n,i}\supset \cdots \supset \ggu_{n,n}\supset \ggu_{n,n+1}:=0
\end{equation}
where $\ggu_{n,i}:=\sum_{j=i}^nP_{j-1}\der_j$ for $i=1, \ldots ,
n$. For all $i<j$, 
\begin{equation}\label{cuni}
[\ggu_{n,i}, \ggu_{n,j}]\subseteq
\begin{cases}
\ggu_{n,i+1}& \text{if }i=j,\\
\ggu_{n,j}& \text{if }i<j.
\end{cases}
\end{equation}
(Proposition 2.1.(2), \cite{Bav-Lie-Un-GEN}) states {\em that (\ref{seruni}) is the derived series for the Lie algebra $\ggu_n$, i.e., $(\ggu_n)_{(i)}=\ggu_{n,i+1}$ for all} $i\geq 0$.

{\bf The group of automorphisms of the Lie algebra $\ggu_n$}.
In \cite{Bav-Lie-Un-AUT}, the group of automorphisms $G_n$ of the Lie algebra $\ggu_n$ of unitriangular polynomial derivations is found ($n\geq 2$), it is isomorphic to an iterated semi-direct  product (Theorem 5.3,  \cite{Bav-Lie-Un-AUT}),
 $$\mT^n\ltimes (\UAut_K(P_n)_n\rtimes( \mF_n' \times  \mE_n ))  $$
 where $\mT^n$ is an algebraic  $n$-dimensional torus,   $\UAut_K(P_n)_n$ is an explicit factor group of the group $\UAut_K(P_n)$ of unitriangular polynomial automorphisms, $\mF_n'$ and $\mE_n$ are explicit groups that are isomorphic respectively to the groups $\mI$ and $\mJ^{n-2}$ where
    $\mI := (1+t^2K[[t]], \cdot )\simeq K^{\N}$ and
  $\mJ := (tK[[t]], +)\simeq K^\N$. It is shown that the {\em adjoint group} of automorphisms $\CA (\ggu_n)$  of the Lie algebra $\ggu_n $ is equal to the group $\UAut_K(P_n)_n$ (Theorem 7.1,  \cite{Bav-Lie-Un-AUT}). Recall that the {\em adjoint group} $\CA (\CG )$ of a Lie algebra $\CG$ is generated by the elements $ e^{ \ad (g)}:=\sum_{i\geq 0}\frac{\ad (g)^i}{i!}\in \Aut_K(\CG )$ where $g$ runs through all the locally nilpotent elements of the Lie algebra $\CG$ (an element $g$ is a {\em locally nilpotent element} if the inner derivation $\ad (g):= [g, \cdot ]$ of the Lie algebra $\CG$ is a locally nilpotent derivation). The group $G_n$ contains the semi-direct product $\mT^n\ltimes \CT_n$ where
  $$\CT_n:=\{ \s \in \Aut_K(P_n)\, | \, \s (x_1)=x_1, \s(x_i) = x_i+a_i\; {\rm where}\;\; a_i\in (x_1, \ldots , x_{i-1}), i=2, \ldots , n\}$$ where $ (x_1, \ldots , x_{i-1})$ is the maximal ideal of the polynomial algebra $P_{i-1}:=K[x_1, \ldots , x_{i-1}]$ generated by the elements $x_1, \ldots , x_{i-1}$.


{\bf Proof of Theorem \ref{10Mar12}}. Let $\v : \ggu_n\ra \ggu_n$ be a monomorphism of the Lie algebra $\ggu_n$. By (Proposition 2.1.(2), \cite{Bav-Lie-Un-GEN}), $(\ggu_n)_{(i)}=\ggu_{n,i+1}$ for all $i$. So, the descending chain of ideals (\ref{seruni}) is the derived series for the Lie algebra $\ggu_n$ of length $l(\ggu_n)=n$ (by definition, this is the number of nonzero terms in the derived series). Clearly, $l(\ggu_{n,2})=n-1$ and
$$l(\v (\ggu_n))=l(\ggu_n) =n$$  $(\v (\ggu_n) \simeq \ggu_n$). It follows that $$\v (\ggu_n) \not\subseteq \ggu_{n,2}$$ since otherwise we would have $n=l(\v (\ggu_n))\leq l(\ggu_{n,2})=n-1$, a contradiction. This means that $\der_1':=\v (\der_1)=\l_1 \der_1+u_1$ for some $\l_1\in K^*$ and $u_i\in \ggu_{n,2}$. We use induction on $i$ to show that
\begin{equation}\label{diui}
\der_i':=\v (\der_i)=\l_i\der_i+u_i,\;\; i=1, \ldots , n,
\end{equation}
for some elements $\l_i\in K^*$ and  $u_i\in \ggu_{n,i+1}$. In particular, $\der_n = \l_n\der_n$. The initial step, $i=1$, has already been established. Suppose that $i\geq 2$ and that (\ref{diui}) holds for all numbers $i'<i$. Since $\v ((\ggu_n)_{(j)})\subseteq (\ggu_n)_{(j)}$ for all $j\geq 1$, we have the inclusion $\v (\ggu_{n,i}) = \v ((\ggu_n)_{(i-1)})\subseteq (\ggu_n)_{(i-1)}=\ggu_{n,i}$ which implies that $\der_i'=\l_i\der_i+u_i$ for some elements $\l_i\in P_{i-1}$ and $u_i\in \ggu_{n,i+1}$. It remains to show that $\l_i\in K^*$. This fact follows from the commutation relations $[\der_j', \der_i']=0$ for $j=1, \ldots , i-1$ ($0=\v ([\der_j, \der_i])=[\der_j', \der_i']$). In more detail, for $j=i-1$,
$$0=[\der_{i-1}', \der_i']=[\l_{i-1}\der_{i-1}+u_{i-1}, \l_i\der_i+u_i]=\l_{i-1}\der_{i-1}(\l_i) \der_i+v_{i-1}$$ for some element $v_{i-1}\in \ggu_{n,i+1}$. Therefore, $\der_{i-1}(\l_i) =0$, i.e.,  $\l_i\in P_{i-2}$. Now, we use a second downward induction on $j$ starting on $j=i-1$ to show that
\begin{equation}\label{diui1}
\l_i\in P_j \;\; {\rm for \; all}\;\; j=1, \ldots , i-1.
\end{equation}
The initial step, $j=i-1$, has been just proved. Suppose that (\ref{diui1}) is true for  all $j=k, \ldots , i-1$. In particular, $\l_i \in P_k=K[x_1, \ldots , x_k]$. We have to show that $\l_i\in P_{k-1}$. For, we use the equality $[\der_k', \der_i']=0$:
$$ 0=[\l_k\der_k+u_k , \l_i\der_i +u_i]=\l_k\der_k(\l_i) \der_i +v_k$$ for some element $v_k\in \ggu_{n,i+1}$ $([u_k,\l_i\der_i]\in \ggu_{n, i+1}$ since $\l_i\in P_k$ and $[\oplus_{k+1\leq j\leq i}P_{j-1}\der_j, \l_i\der_i]=0$). Therefore, $\der_k(\l_i)=0$, i.e.,  $\l_i\in P_{k-1}$. By induction on $j$,  (\ref{diui1}) holds. In particular, for $j=1$:
 $\l_i\in P_{1-1}=P_0=K$. We have to show that $\l_i\neq 0$. Notice that $\ggu_{n,i}=\oplus_{j=i}^nP_{j-1}\der_j$, $l(\ggu_{n,i}) = n-i+1$ and $ \v ( \ggu_{n,i}) \subseteq \ggu_{n,i}$. The monomorphism $\v$ respects the Lie subalgebra $\CG = K\der_i+\ggu_{n,i+1}$ of the Lie algebra $\ggu_n$, i.e., $\v (\CG ) \subseteq \CG$. The inclusion of Lie algebras $\CG \subseteq \ggu_{n,i}$ yields the inequality $l(\CG ) \leq l(\ggu_{n,i})=n-i+1$ (the equality follows from the fact that
  $(\ggu_n)_{(j)}=\ggu_{n,j+1}$ for all $j\geq 0$). The vector space
  $$ \CH = K\der_i+K[x_i]\der_{i+1}+K[{x_i}, x_{i+1}]\der_{i+2}+\cdots + K[x_i, \ldots , x_{n-1}]\der_n$$
is a Lie subalgebra of $\CG$ which is isomorphic to the Lie algebra $\ggu_{n-i+1}$. Therefore, $ l(\CH ) = l(\ggu_{n-i+1}) = n-i+1$. The inclusion of Lie algebra $\CH \subseteq \CG$ yields the inequality $n-i+1= l(\CH ) \leq l(\CG )$. Therefore, $l(\CG ) = n-i+1$.

  Suppose that $\l_i=0$, we seek a contradiction. In that case, $\v (\CG ) \subseteq \ggu_{n,i+1}$ and so
$$ n-i+1=l(\CG )=l(\v (\CG )) \leq l(\ggu_{n,i+1})=n-i,$$
 a contradiction.

  Therefore, (\ref{diui}) holds. By (Theorem 3.6.(2), \cite{Bav-Lie-Un-AUT}), there exists a unique  automorphism $\s \in \mT^n \ltimes \CT_n\subseteq G_n$ such that $\s (\der_i) = \der_i'$ for $i=1, \ldots , n$. By replacing the monomorphism $\v$ by the monomorphism $\s^{-1}\v$, without loss of generality we can assume that $$\der_i'=\der_i\;\; {\rm  for \; all }\;\;  i=1, \ldots , n.$$ The vector space $\ggu_n$ is the union $\cup_{i\geq 0} N_i$ of vector subspaces $N_i:= \{ u\in \ggu_n\, | \, \d_j^{i+1}(u)=0, \; j=1, \ldots , n-1\}$ where $\d_j=\ad (\der_j)$. Clearly, $N_i= \oplus_{j=1}^n N_i\cap P_{j-1}\der_j$ and $ N_i\cap P_{j-1}\der_j=\oplus\{ Kx^\alpha \, | \, \alpha = (\alpha_1, \ldots , \alpha_{j-1})\in \N^{j-1}, \alpha_k\leq i$ for $k=1, \ldots , j-1\}$. In particular, $$\dim_K(N_i)<\infty\;\; {\rm  for\;  all}\;\;  i\geq 0.$$ By the very definition of the vector spaces $N_i$ and the fact that $\v (\der_i) = \der_i$ for $i=1, \ldots , n-1$, $$\v (N_i) \subseteq N_i\;\; {\rm  for \;  all}\;\;  i\geq 0.$$ Since the linear map $\v$ is an injection and the vector spaces $N_i$ are finite dimensional, we have $\v (N_i) = N_i$ for all $i\geq 0$, i.e.,  $\v$ is a bijection. $\Box $


$${\bf Acknowledgements}$$

 The work is partly supported by  the Royal Society  and EPSRC.

\small{

Department of Pure Mathematics

University of Sheffield

Hicks Building

Sheffield S3 7RH

UK

email: v.bavula@sheffield.ac.uk}


\begin{thebibliography}{99}
\bibitem{Bav-cdixintdif}  V. V. Bavula,  An analogue of the Conjecture of Dixmier is true for the
algebra
 of polynomial integro-differential operators,  Arxiv:math.RA: 1011.3009.
\end{thebibliography}
\end{document}